 \magnification=1200
\input amssym.def
\input amssym.tex
 \font\newrm =cmr10 at 24pt
\def\bul{\raise .9pt\hbox{\newrm .\kern-.105em } }
 \def\fr{\frak}

 \baselineskip 20pt

 \def\h{\hbox{ }}

 \def\n{{\fr n}}
 \def\a{{\fr a}}

 \def\ss{{\fr s}}
 
 \def\b{{\fr b}}
 \def\cc{{\fr c}}
 \def\hh{{\fr h}}

 \def\g{{\fr g}}

 \def\<{\le}
 \def\>{\ge}

 \def\s{{\h\subset\h}}

 \def\vs{\vskip }

 \def\mapright#1
  {\smash{\mathop
  {\longrightarrow}
  \limits^{#1}}}

 \def\kk#1{{\kern .4 em} #1}
 \def\vs{\vskip 1pc}

\font\twelvebf  = cmb10 scaled 1400
\hsize = 31pc
\vsize = 45pc
\overfullrule = 0pt
\font\srm   =cmr10    scaled 700  

\rm

\centerline{\twelvebf On the algebraic set of singular elements}
\centerline{\twelvebf in a
complex simple Lie algebra}
\vs

\centerline{BERTRAM KOSTANT  and NOLAN WALLACH}
\vs
\noindent {\bf Abstract.} Let $G$ be a complex simple Lie group and let $\g =
\hbox{\rm Lie}\,G$.  Let $S(\g)$ be the $G$-module of polynomial functions on $\g$ and let
$\hbox{\rm Sing}\,\g$ be the closed algebraic cone of singular elements in $\g$. Let
${\cal L}\s S(\g)$ be the (graded) ideal defining $\hbox{\rm Sing}\,\g$ and let $2r$ be the
dimension of a $G$-orbit of a regular element in $\g$. Then ${\cal L}^k = 0$ for any
$k<r$. On the other hand, there exists a remarkable $G$-module $M\s {\cal L}^r$ which
already defines $\hbox{\rm Sing}\,\g$. The main results of this paper are a determination of
the structure of $M$.

\vskip 1.5pc \centerline{\bf 0. Introduction} \vskip.5pc
{\bf 0.1.} Let $G$ be a complex simple Lie group and let $\g =
\hbox{\rm Lie}\,G$. Let $\ell = \hbox{\rm rank}\,\g$. Then in superscript centralizer
notation one has $\hbox{\rm dim}\,\g^x \geq \ell$ for any $x\in \g$. An
element $x\in \g$ is called regular (resp. singular) if $\hbox{\rm dim}\,\g^x
= \ell$ (resp. $ > \ell$). Let $\hbox{\rm Reg}\,\g$ be the set of all regular
elements in $\g$ and let $\hbox{\rm Sing}\,\g$, its complement in $\g$, be
the set of all singular elements in $\g$. Then one knows that
$\hbox{\rm Reg}\,\g$ is a nonempty Zariski open subset of $\g$ and hence
$\hbox{\rm Sing}\,\g$ is a closed proper algebraic subset of $\g$.

Let $S(\g)$ (resp. $\wedge\,\g$) be the symmetric (resp. exterior)
algebra over $\g$. Both algebras are graded and are $G$-modules by
extension of the adjoint representation. Let ${\cal B}$ be the
natural extension of the Killing form to $S(\g)$ and $\wedge\,\g$.
The inner product it induces on $u$ and $v$ in either $S(\g)$ or
$\wedge\,\g$ is denoted by $(u,v)$. The use of ${\cal B}$ permits
an identification of $S(\g)$ with the algebra of polynomial
functions on $\g$. Since $\hbox{\rm Sing}\,\g$ is clearly a cone the ideal,
${\cal L}$, of all $f\in S(\g)$ which vanish on $\hbox{\rm Sing}\,\g$ is graded.
Let $n = \hbox{\rm dim}\, \g$ and let $r = (n-\ell)/2$. One knows that $n-\ell$ is
even so that $r\in {\Bbb Z}_+$. It is easy to show that $${\cal L}^k
= 0,\,\,\,\hbox{for all $k<r$.}\eqno (0.1)$$
 The purpose of this
paper is to define and study a rather remarkable $G$-submodule $$M\s
{\cal L}^r\eqno (0.2)$$ which in fact defines $\hbox{\rm Sing}\,\g$. That is,
if
$x\in \g$, then $$x\in\,\hbox{\rm Sing}\,\g \iff\,\,f(x)= 0,\,\,\,\forall f\in
M\eqno (0.3)$$ \vskip .5pc
{\bf 0.2.} We will now give a definition of $M$. The use of ${\cal
B}$ permits an identification of $\wedge\,\g$ with the underlying
space of the cochain complex defining the cohomology of $\g$. The
coboundary operator is denoted here by $d$ (and $\delta$ in [Kz])
is a (super) derivation of degree 1 of $\wedge\,\g$ so that $dx\in
\wedge^2\g$ for any $x\in \g$. Since $\wedge^{\hbox{\srm even}}\g$ is a
commutative algebra there exists a homomorphism $$\gamma:S(\g)\to
\wedge^{\hbox{\srm even}}\g$$ where for $x\in \g$, $\gamma(x) = -dx$. One
readily has that $$S^k(\g)\s \hbox{\rm Ker}\,\gamma,\,\,\hbox{for all $k>r$.}
\eqno (0.4)$$ Let $\gamma_r = \gamma|S^r(\g)$ so that
$$\gamma_r:S^r(\g)\to
\wedge^{2r}\g.\eqno (0.5)$$ If
$x\in \g$, one readily has
$$x^r\in \hbox{\rm Ker}\,\gamma_r \,\,\iff x\in \hbox{\rm Sing}\,\g. \eqno
(0.6)$$\vskip.5pc
Let $\Gamma$ be the transpose of $\gamma_r$ so
that one has a $G$-map $$\Gamma:\wedge^{2r}\to S^r(\g). \eqno (0.7)$$
By definition $$M = \hbox{\rm Im}\,\,\Gamma.\eqno (0.8)$$	\vskip .5pc

{\bf 0.3.} Let $J= S(\g)^G$ so that (Chevalley) $J$ is a polynomial
ring ${\Bbb C}[p_1,\ldots,p_{\ell}]$ where the invariants $p_j$ can
be chosen to be homogeneous. In fact if
$m_j,\,j=1,\dots,\ell$, are the exponents of $\g$ we can take
$\hbox{\rm deg}\,p_j = m_j +1$. For any linearly independent
$u_1,\ldots,u_{\ell}\in \g$, let $$\psi(u_1,\ldots,u_{\ell}) =
\hbox{\rm det}\,\partial_{u_i}\,p_j \eqno (0.9)$$ where, if $v\in \g$,
$\partial_v$ is the operator of partial derivative by $v$ in $S(\g)$.
One has $$\psi(u_1,\ldots,u_{\ell}) \in S^r(\g)\eqno (0.10)$$ since,
as one knows, $\sum_{i=1}^{\ell} m_i = r$.

Let $\Sigma_{2r}$ be the permutation group of $\{1,\ldots,2r\}$ and let
$\Pi_r\s \Sigma_{2r}$ be a subset (of cardinality $(2r-1)(2r-3)\cdots
1$) with the property that $sg\,\nu = 1$ for all $\nu\in \Pi_r$ and such
that, as unordered, $$\{(\nu(1),\nu(2)),\ldots,
(\nu(2r-1),\nu(2r))\}\mid \nu\in \Pi_r\}$$ is the set of all partitions
of $\{1,\ldots,2r\}$ into a union of $r$ subsets each of which has two
elements. The following is one of our main theorems. Even more than
explicitly determining $\psi(u_1,\ldots,u_{\ell})$ one has

  \vs {\bf
Theorem 0.1.} {\it Let $u_1,\ldots,u_{\ell}$ be any $\ell$ linearly
independent elements in $\g$ and let $w_1,\dots,w_{2r}$ be a basis of
the ${\cal B}$-orthogonal subspace to the span of the $u_i$. Then there
exists some fixed $\kappa\in
\hbox{\rm \Bbb C}^{\times}$ such that, for all $x\in \g$,
$$\sum_{\nu\in \Pi_r}([w_{\nu(1)},w_{\nu(2)}],x)\cdots
([w_{\nu(2r-1)},w_{\nu(2r)}],x)= \kappa\,\,
\psi(u_1,\ldots,u_{\ell})(x).\eqno (0.11)$$ Moreover
$\psi(u_1,\ldots,u_{\ell}) \in M$. In fact the left side of (0.11) is
just $\Gamma(w_1\wedge\cdots\wedge w_{2r})(x)$. In addition $M$ is the
span of $\psi(u_1,\ldots,u_{\ell})$, over all $\{u_1,\ldots,u_{\ell}\}$,
taken from the ${n\choose l}$ subsets of $\ell$-elements in any given
basis of $\g$.}\vs

We now deal with the $G$-module structure of $M$. For
any subspace $\ss$ of $\g$, say of dimension $k$, let $[\ss] = \Bbb C
v_1\wedge \cdots \wedge v_k\s \wedge^k\g$ where the $v_i$ are a basis of
$\ss$. Let
$\hh$ be a Cartan subalgebra of $\g$ and let $\Delta$ be the set of roots
for the pair $(\hh,\g)$. For any $\varphi\in \Delta$ let $e_{\varphi}\in
\g$ be a corresponding root vector. Let $\Delta_+\s\Delta$ be a choice of
a set of positive roots and let $\b$ be the Borel subalgebra spanned by
$\hh$ and all $e_{\varphi}$ for $\varphi\in \Delta_+$. For any subset
$\Phi\s \Delta$ let $\a_{\Phi}\s \g$ be the span of $e_{\varphi}$ for
$\varphi\in \Phi$. Also let $\langle \Phi\rangle = \sum_{\varphi\in \Phi}
\varphi$ so that $$[\a_{\Phi}]\,\,\hbox{is an $\hh$-weight space for the
$\hh$-weight $\langle \Phi\rangle$.} \eqno (0.12)$$\vskip .5pc

 A subset
$\Phi\in
\Delta_+$ will be said to be an ideal in $\Delta_+$ if $\a_{\Phi}$ is an
ideal of $\b$. In such a case, if $\hbox{\rm card}\,\Phi = k$, then the span
$V_{\Phi}$ of $G\cdot [a_{\Phi}]$ is an irreducible $G$-submodule of
$\wedge^k\g$ having $[a_{\Phi}]$ as highest weight space and $\langle
\Phi\rangle $ as highest weight. Let ${\cal I}$ be the set of all ideals
$\Phi$ in $\Delta_+$ of cardinality $\ell$. It is shown in [KW]
that all ideals in $\b$ of dimension $\ell$ are abelian and hence are
of the form $\a_{\Phi}$  for a unique $\Phi\in {\cal I}$. Specializing $k$ in
[K3] to $\ell$ one has that, by definition, $A_{\ell}\s
\wedge^{\ell}\g$ is the span of $[\ss]$ over all abelian subalgebras $\ss\s \g$ of
dimension
$\ell$. Using results in [K3] and that in [KW] above, one also has that
$A_{\ell}$ is a multiplicity one $G$-module with the complete reduction
$$A_{\ell} =
\oplus_{\Phi\in {\cal I}} V_{\Phi}\eqno (0.13)$$ so that there are exactly
$\hbox{\rm card}\,{\cal I}$ irreducible components. In addition it has been shown in [K3]
that $\ell$ is the maximal eigenvalue of the (${\cal B}$ normalized)
Casimir operator,
$\hbox{\rm Cas}$, in
$\wedge^{\ell}\g$ and $A_{\ell}$ is the corresponding eigenspace. In the
present paper the $G$-module structure of $M$ is given in \vs {\bf Theorem
0.2.} {\it As
$G$-modules one has an equivalence $$M\cong A_{\ell}\eqno (0.14)$$ so
that $M$ is a multiplicity one module with $\hbox{\rm card}\,{\cal I}$ irreducible
components. Morever the components can be parameterized by ${\cal I}$ in
such a way that the component corresponding to $\Phi\in {\cal I}$ has
highest weight $\langle
\Phi\rangle$. In addition $\hbox{\rm Cas}$ takes the value $\ell$ on each and every
irreducible component of $M$. }\eject

\centerline{\bf 1. Preliminaries}\vs {\bf 1.1.} Let $\g$ be
a complex semisimple Lie algebra and let $G$ be a Lie group such
that $\g = \hbox{\rm Lie}\,G$. Let $\hh\s \g$ be a Cartan subalgebra of $\g$
and let $\ell$ be the rank of $\g$ so that $\ell = \hbox{\rm dim}\,\hh$. Let
$\Delta$ be the set of roots for the pair $(\hh,\g)$ and let
$\Delta_+\s\Delta$ be a choice of a set of positive roots. Let $r =
\hbox{\rm card}\,\Delta_+$ so that $$n = \ell + 2\,r\eqno (1.1)$$ where we
let $n= \hbox{\rm dim}\,\g$. Let ${\cal B}$ be Killing form $(x,y)$ on $\g$.
For notational economy we identify $\g$ with its dual $\g^*$ using
${\cal B}$. The bilinear form ${\cal B}$ extends to an inner product
$(p,q)$, still denoted by ${\cal B}$, on the two graded algebras,
the symmetric algebra
$S(\g)$ of $\g$ and the exterior algebra $\wedge\,\g$ of $\g$. If
$x_i,y_j\in
\g,\, i=1,\ldots,k,\,j=1,\ldots,m$, then the product of $x_i$ is
orthogonal to the product of $y_j$ in both $S(\g)$ and $\wedge\,\g$
if $k\neq m$, whereas if $k=m$, $$\eqalign{(x_1\cdots x_k,y_1\cdots
y_k)&= \sum_{\sigma\in \Sigma_k} (x_1,y_{\sigma(1)})\cdots
(x_k,y_{\sigma(k)})\,\,\,\,\,\hbox{in}\,\,S(\g)\cr
(x_1\wedge\cdots \wedge x_k,y_1\wedge \cdots
\wedge y_k)&= \sum_{\sigma\in \Sigma_k}
sg(\sigma) (x_1,y_{\sigma(1)})\cdots
(x_k,y_{\sigma(k)})\,\,\,\,\,\hbox{in}\,\,\wedge\g.\cr}\eqno
(1.2)$$ Here $\Sigma_k$ is the permutation group on
$\{1,\ldots,k\}$ and $sg$ abbreviates the signum character on $\Sigma_k$.

The identification of $\g$ with its dual has the effect of
identifying $S(\g)$ with the algebra of polynomial functions $f(y)$ on
$\g$. Thus if $x,y\in \g$, then $x(y) = (x,y)$ and if $x_i\in
\g,\,i=1,\ldots,k$, then $$\eqalign{(x_1\cdots
x_k)(y)&=\prod_{i=1}^k (x_i,y)\cr &= (x_1\cdots x_k,{1\over
k!}\,y^k).\cr}\eqno (1.3)$$

\vskip .5pc The identification of
 $\g$ with its dual also has the effect of
identifying the (supercommutative) algebra $\wedge\,\g$ with the
underlying space of the standard cochain complex defining the
cohomology of $\g$. Let $d$ be the (super) derivation of degree 1
of $\wedge\g$, defined by putting $$d = {1\over 2}\,\sum_{i=1}^n
\varepsilon(w_i)\theta (z_i).\eqno (1.4)$$ Here $\varepsilon(u)$, for
any $u\in \wedge \g$, is left exterior multiplication by $u$ so
that $\varepsilon(u)\,v = u\wedge v$ for any $v\in \wedge\,\g$.
Also $w_i, i=1,\ldots,n$, is any basis of $\g$ and
$z_i\in \g,\,i=1,\ldots,n$, is the ${\cal B}$ dual basis.
$\theta(x)$, for
$x\in\g$, is the derivation of $\wedge\,g$, of degree 0, defined
so that $\theta(x)y = [x,y]$ for any $y\in \g$. One readily notes
that (1.4) is independent of the choice of the basis $w_i$. Thus
if $x\in \g$, then $dx\in \wedge^2\g$ is given by $$dx = {1\over 2}
\sum_{i=1}^n w_i\wedge [z_i,x]. \eqno (1.5)$$
\vskip .5pc
Any element $\omega\in \wedge^2\g$ defines an alternating
bilinear form on $\g$. Its value $\omega(y,z)$ on $y,z\in \g$ may
be given in terms of ${\cal B}$ by $$\omega(y,z) = (\omega, y\wedge
z).\eqno (1.6)$$  The rank of $\omega$ is necessarily even. In
fact if $\hbox{\rm rank}\,\omega = 2k$, then there exist $2k$ linearly
independent elements $v_i\in \g,\,i=1,\ldots,2k$, such that $$\omega
= v_1\wedge v_2 + \cdots + v_{2k-1}\wedge v_{2k}.\eqno (1.7)$$ The
radical of $\omega$, denoted by $\hbox{\rm Rad}\,\omega$, is the space of all
$y\in \g$ such that $\omega(y,z) = 0$ for all $z\in \g$. For
$u\in \wedge\,\g$, let $\iota(u)$ be the transpose of
$\varepsilon(u)$ with respect to ${\cal B}$ on $\wedge\,\g$. If
$u=y\in\g$, then one knows that $\iota(y)$ is the (super) derivation
of degree minus 1 defined so that if $z\in \g$, then $\iota(y)z = (y,z)$.
(See p.~8 in [Kz]). From (1.6) one has $$\hbox{\rm Rad}\,\omega = \{y\in
\g\mid \iota(y)\omega = 0\}.\eqno (1.8)$$ If $\ss$ is any
subspace of $\g$, let $\ss^{\perp}$ be the ${\cal B}$ orthogonal
subspace to $\ss$. From (1.7) one then has that
$$\{v_i\},\,i=1,\ldots,2k, \,\,\,\hbox{is a basis of
$\hbox{\rm Rad}\,\omega^{\perp}$}. \eqno (1.9)$$ If $\ss\s \g$ is any subspace,
say of dimension $m$, let $[\ss]\in \wedge^m\g$ be the $\Bbb C$
span of the decomposable element $u_1\wedge \cdots\wedge u_m$
where $\{u_i,\,i=1,\ldots,m\}$ is a basis of $\ss$. One notes
that if $\omega\in \wedge^2\g$ is given as in (1.7), then $$\omega^k
= k!\,\,\,v_1\wedge\cdots \wedge v_{2k}\eqno (1.10)$$ so that
$$\omega^j\neq 0 \iff j\leq k\,\,\,\hbox{and}\,\, \omega^k \in
[\hbox{\rm Rad}\,\omega^{\perp}].\eqno (1.11)$$

\vskip .5pc Let
$\{w_j,\,j=1,\ldots,n\}$ be a ${\cal B}$ orthonormal basis of $\g$.
Put $\mu = w_1\wedge\cdots\wedge w_n$ so that $$(\mu,\mu) = 1\eqno
(1.12)$$ so that $\mu$ is unique up to sign and $\wedge^n\g =\Bbb
C\mu$. For any $v\in \wedge\,\g$ let $v^* = \iota(v)\mu$. We recall
the more or less well known.

 \vs {\bf Proposition 1.1.} {\it If
$\ss\s \g$ is any subspace and $0\neq u\in [\ss]$, then $$0\neq
u^*\in [\ss^{\perp}].\eqno (1.13)$$ Moreover if $s,t\in \wedge\,\g$,
one has $$(s,t) = (s^*,t^*).\eqno (1.14)$$ }\vskip.5pc
{\bf Proof.} Let
$\{y_i,\,\,i=1,\ldots,m\}$ be a basis of $\ss$ chosen so that $u
= y_1 \wedge \cdots \wedge y_m$ and let
$\{z_j,\,\,j=1,\ldots,n-m\}$ be a basis of $\ss^{\perp}$. Then if
$y'_k,\,\,k=1,\ldots,m$, are chosen in $\g$ such that $(y_i,y'_k) =
\delta_{ik}$, it is immediate that the $y_k'$ together with the
$z_j$ form a basis of $\g$ so that for some $\lambda\in \Bbb
C^{\times}$ one has $$\lambda y'_1\wedge \cdots\wedge y'_m\wedge
z_1\wedge \cdots\wedge z_{n-m} = \mu.\eqno (1.15)$$ But since
interior product is the transpose of exterior product one has
$$\iota(q)\,\iota(p) = \iota(p\wedge q)\eqno(1.16)$$ for any
$p,q\in \wedge^\g$. Thus by (1.15) one has $$u^* = \lambda
z_1\wedge \cdots\wedge z_{n-m}$$ establishing (1.13). To prove
(1.14) it suffices by linearity to assume that both $s$ and $t$ are
decomposable of some degree $m$. Thus we can assume $s = y_1\wedge
\cdots \wedge y_m$ and $t= z_1\wedge \cdots \wedge z_m$ for
$y_i,z_j\in \g$. But now, as one knows, and readily establishes,
$$\varepsilon(y)\,\iota (z) + \iota (z) \varepsilon (y) = (y,z)
\hbox{\rm Id}_{\g}\eqno (1.17)$$ for $y,z\in \g$. Thus
$$\eqalign{(s^*,t^*)&=(\iota(s)\mu,\iota(t)\mu)\cr &=(\mu,
\varepsilon(s)\iota(t)\mu)\cr}.\eqno (1.18)$$ But then using
(1.17) and the fact that $\varepsilon(y)\mu= 0$ for any $y\in \g$, one
has $$(\mu,\varepsilon(s)\iota(t)\mu)
= \sum_{j=0}^{m-1}
(-1)^{j}(y_m,z_{m-j})(\mu,\varepsilon(y_1)\cdots\varepsilon(y_{m-1})\,
\iota(z_m)\cdots {\widehat {\iota(z_{m-j})}}\cdots
\iota(z_1)\mu).$$ But then by induction and the expansion of the
determinant defined by the last row one has
$$\eqalign{(\mu,\varepsilon(s)\iota(t)\mu) &=
det\,(y_i,z_j)(\mu,\mu)\cr &= (s,t)\cr}$$ proving (1.14). QED

\vs
{\bf 1.2.} The algebra $S(\g)$ is a $G$-module extending the adjoint
representation. Let $J = S(\g)^{G}$ be the subalgebra of $\g$-invariants.
Let $H\s S(\g)$  be the graded $\g$-submodule of harmonic elements in
$S(\g)$ (See  \S 1.4 in [K2] for definitions). Then one knows
$$S(\g) = J\otimes H.\eqno (1.19)$$
See  (1.4.3) in [K2].

Let $r$ be as in (1.1). For the convenience of the reader we repeat a
paragraph in \S 1.2 of [K4]. Let
$\Sigma_{2r,2}$ be the subgroup of all
$\sigma\in \Sigma_{2r}$ such that $\sigma$ permutes the set of unordered
pairs $\{(1,2),(3,4),\ldots,(2r-1,2r)\}$. It is clear that
$\Sigma_{2r,2}$ has order $r!\,2^r$. Now let $\Pi_r$ be a cross-section
of the set of left cosets of $\Sigma_{2r,2}$ in $\Sigma_{2r}$. Thus one
has a disjoint $$\Sigma_{2r} = \bigcup_{\nu\in \Pi_r} \nu
\,\Sigma_{2r,2}.\eqno (1.20)$$ One notes that the cardinality of $\Pi_r$
is $(2r-1)(2r-3)\cdots 1$ (the index of $\Sigma_{2r,2}$ in
$\Sigma_{2r}$) and the correspondence $$\nu\mapsto
((\nu(1),\nu(2)),(\nu(3),\nu(4)),\ldots,(\nu(2r-1),\nu(2r))\eqno
(1.21)$$ sets up a bijection of $\Pi_r$ with the set of all partitions
of $(1,2,\ldots,2r)$ into a union of subsets, each of which has two
elements. Furthermore, since the signum character restricted to
$\Sigma_{2r,2}$ is nontrivial we may choose $\Pi_r$ so that $$sg(\nu) =
1$$ for all $\nu\in \Pi_r$.

In [K4] we defined a map $\Gamma:\wedge^{2r}\g \to S(\g)$;  (Its
significance will  become apparent later). Here, using Proposition 1.2
in [K4] we will give a simpler definition of
$\Gamma$. By Proposition 1.2 in [K4] one has \vs {\bf Proposition
1.2.} {\it There exists a map $$\Gamma:\wedge^{2r}\g \to S^r(\g)\eqno
(1.21a)$$ such that for any $w_i\in\g,\,i=1,\ldots,2r$, one has
$$\Gamma(w_1\wedge \cdots \wedge w_{2r}) = \sum_{\nu\in
\Pi_r}[w_{\nu(1)},w_{\nu(2)}]\cdots [w_{\nu(2r-1)},w_{\nu(2r)}].\eqno
(1.22)$$}\vs As a polynomial function of degree $r$ on $\g$,
one notes that $$\Gamma(w_1\wedge \cdots\wedge w_{2r})(x) =
\sum_{\nu\in
\Pi_r}([w_{\nu(1)},w_{\nu(2)}],x)\cdots
([w_{\nu(2r-1)},w_{\nu(2r)}],x).\eqno (1.23)$$ This clear from (1.1.7)
in [K4] and (1.3) here.

The algebra $\wedge\g$ is a natural $G$-module by extension of the
adjoint representation. It is clear that $\Gamma$ is a $G$-map. Let $M\s
S^r(\g)$ be the image of $\Gamma$. The following is proved as Corollary
3.3 in [K4]. \vs {\bf Theorem 1.3.} {\it One has
$M\s H^r$ so that $M$ is a $G$-module of harmonic polynomials of
degree $r$ on
$\g$.}\vs Giving properties of $M$ and determining its
rather striking $\g$-module structure is the main goal of this paper.

For any $y\in \g$ one has the familiar supercommutation formula
$\iota(y)d + d \iota(y) = \theta(y)$. See e.g., (92) in [K5]. Now let
$x,y\in \g$. Since $d \iota(y) (x)= 0$ one has $\iota(y)dx = [y,x]$.
Thus, by (1.8), using superscript notation for centralizers one has $$\hbox{\rm Rad}\,dx
= \g^x. \eqno (1.24)$$ Clearly $[x,\g]$ is the ${\cal B}$ orthogonal
subspace in $\g$ to $\g^x$ so that $$[x,\g] = (\hbox{\rm Rad}\,dx)^{\perp}\eqno (1.25)$$
for any $x\in \g$.

For any $x\in \g$ one knows $\hbox{\rm dim}\,\g^x \geq \ell$. Recall that an element $x\in \g$ is
called regular if $\hbox{\rm dim}\,\g^x = \ell$. The set $\hbox{\rm Reg}\,\g$ of regular
 elements
is nonempty and Zariski open. Its complement, $\hbox{\rm Sing}\,\g$, is the Zariski closed
set of singular elements. One notes, by (1.11), that $$\hbox{\rm Sing}\,\g = \{x\in
\g\mid (dx)^r = 0\}.\eqno (1.26)$$\vskip .5pc
Now
$\wedge^{\hbox{\srm even}}\g$ is a commutative algebra and hence there exists a
homomorphism $$\gamma:S(g)\to
\wedge^{\hbox{\srm even}}\g\eqno (1.27)$$ such that for $x\in \g$, $$\gamma(x) = -dx.$$ Let
$\gamma_r$ be the restriction of $\gamma$ to $S^r(\g)$. The following result,
established as Theorem 1.4 in [K4], asserts that $\Gamma$ is the transpose of
$\gamma_r$.

\vs {\bf Theorem 1.4.} {\it Let $y_1,\ldots, y_r\in \g$ and let
$\zeta\in
\wedge^{2r}(\g)$. Then $$(y_1\cdots y_r,\Gamma(\zeta)) =
(-1)^r(dy_1\wedge\cdots \wedge dy_r, \zeta).\eqno (1.28)$$}\vskip .5pc Now one knows
that $S^r(\g)$ is (polarization) spanned by all powers $x^r$ for $x\in \g$.
Using (1.3), (1.26) and Theorem 1.4 we recover Proposition 3.2 in [K4]. The
key point is that $M$ defines the variety $\hbox{\rm Sing}\,\g$. \vs {\bf Theorem
1.5.}  {\it Let $x\in \g$ and $\zeta\in \wedge^{2r}\g$. Then
$$\Gamma(\zeta)(x) = {(-1)^r\over r!} ((dx)^r,\zeta).\eqno (1.29)$$ In
particular $$f(x)= 0,\,\,\forall\,f\in M \iff x\in \hbox{\rm Sing}(\g).\eqno (1.30)$$}\vskip.5pc

If $\a$ is a Cartan subalgebra of $\g$,  then one knows that $\a\cap \hbox{\rm Sing}\,\g$
is a union of the root hyperplanes in $\a$. Hence as a corollary of
Theorem 1.5 one has \vs {\bf Theorem 1.6.} {\it Let $\a$ be a Cartan
subalgebra of $\g$. Let $\Delta_+(\a)$ be a choice of positive roots for the
pair $(\a,\g)$. Then for any $f\in M$ one has
$$f|\a\in \Bbb C \prod_{\beta\in \Delta_+(\a)} \beta.\eqno (1.31)$$} \vskip.5pc
 Going
to the opposite extreme we recall that a nilpotent element $e$ is called
principal if it is regular. Let $e$ be a principal nilpotent element. Then by
Corollary 5.6 in [K1] there exists a unique nilpotent radical $\n$ of a
Borel subalgebra such that $e\in \n$. Furthermore $\g^e\cap [\n,\n]$ is a
linear hyperplane in $\g^e$ and $\g^e\cap [\n,\n]= \hbox{\rm (Sing}\,\g)\cap\g^e$ by Theorem
5.3 and Theorem 6.7 in [K1]. Thus there exists a nonzero linear functional
$\xi$ on $\g^e$ such that $$\hbox{\rm Ker}\,\xi = \hbox{\rm (Sing}\,\g)\cap \g^e.\eqno (1.32)$$
 This
establishes \vs {\bf Theorem 1.7.} {\it Let $e\in \g$ be principal nilpotent.
Let $f\in M$. Then using the notation of (1.32) one has $$f|\g^e \in \Bbb
C\,\xi^r. \eqno (1.33)$$}

  Since $\hbox{\rm Sing}\,\g$ is
clearly a cone it follows that the ideal ${\cal L}$ of $f\in S(g)$ which vanishes on
$\hbox{\rm Sing}\,\g$ is graded. One of course has that $M\s {\cal L}^r$.
  We now observe that $r$
is the minimal value of $k$ such that ${\cal L}^k\neq 0$\vs

{\bf Proposition 1.8.} {\it Assume that $0\neq f\in {\cal L}^k$.
Then
$k
\geq r$.}\vs {\bf Proof.} Since $f\neq 0$ there clearly exists a
Cartan subalgebra
$\a$ of $\g$ such that $f|\a \neq 0$. But then
using the notation of Theorem 1.6 it follows from the prime
decomposition that
$\beta$ divides $f|\a$ for all $\beta\in \Delta_+(\a)$.
Thus $k\geq r$. QED

\vskip 1.5 pc \centerline{\bf 2. The structure of $M$
in terms of minors and as a $G$-module} \vskip.5pc {\bf 2.1.} For any
$z\in
\g$ let $\partial_z$ be the partial derivative of $S(\g)$ defined by $z$. Let
$W(\g)= S(\g)\otimes \wedge\,\g$ so that $W(\g)$ can be regarded
as the supercommutative algebra of all differential forms on $\g$ with polynomial
coefficients. To avoid confusion with the already defined $d$, let $d_W$ be the
operator of exterior differentiation on $W(\g)$. That is, $d_W$ is a derivation of
degree 1 defined so that if $\{z_i,w_j\},\,\, i,j = 1,\ldots,n$,  are dual ${\cal B}$
bases  of $\g$, then $$d_W(f\otimes u) = \sum_i^n \partial_{z_i}\,f\otimes
\varepsilon(w_i)\,u\eqno (2.1)$$ where $f\in S(\g)$ and $u\in \wedge\,\g$. Of course
$d_W$ is independent of the choice of bases. In particular $d_Wf$ is a differential
form of degree 1 on $\g$.

For any $x\in \g$ one has a homomorphism $$W(\g)\to \wedge
\,\g,\,\,\,\,\,\,\varphi\mapsto
\varphi(x)\eqno (2.2)$$ defined so that if $\varphi= f\otimes u $, using the notation
of (2.1), then $\varphi(x) = f(x)u$. Next one notes that the  $G$-module
structures on $S(\g)$ and $\wedge\,\g$ define, by tensor product, a  $G$-module
structure on $W(\g)$. Clearly $d_W$ is a $G$  map. If
$a\in G$ and
$\varphi\in W(\g)$, the action of $a$ on $\varphi$ will simply be denoted by $a\cdot
\varphi$. If
$x\in \g$ one readily has
$$a\cdot (\varphi(x)) = a\cdot\varphi (a\cdot x).\eqno (2.3)$$\vskip .5pc

 One knows (Chevalley) that $J$ is a polynomial ring $\Bbb C[p_1,\ldots,p_{\ell}]$
where the $p_j$ are homogeneous polynomials. If $d_j= \hbox{\rm deg}\,p_j$, for
$j=1,\ldots,\ell$, and $m_j=d_j - 1$, then the $m_j$ are exponents of $\g$ so that
$$\sum_{j=1}^{\ell} m_j = r. \eqno (2.4)$$  Moreover we can choose the $p_j$
so that $\partial_y p_j \in H$ for any $y\in \g$ (see Theorem 67 in [K5]). In fact,
if $H_{\hbox{\srm ad}}$ is the primary component of $H$ corresponding to the adjoint
 representation,
then the multiplicity of the adjoint representation in $H_{\hbox{\srm ad}}$ is equal to $\ell$
and
 $\tau_j,\,j=1,\ldots,\ell$, is a basis of $\hbox{\rm Hom}_G(\g,H_{\hbox{\srm ad}})$ where
$$\tau_j(y) =
\partial_y p_j \eqno (2.5)$$ for any $y\in \g$. Again see Theorem 67 in [K5]. \vs

{\bf Remark 2.2.} Using the notation of (2.1) note that $$\{w_{i_1}\wedge \cdots
\wedge w_{i_{\ell}}\mid 1\leq i_1< \cdots < i_{\ell}\leq n\}$$ is a basis of
$\wedge^{\ell}\g$. Furthermore  $$\{z_{j_1}\wedge \cdots
\wedge z_{j_{\ell}}\mid 1\leq j_1< \cdots < j_{\ell}\leq n\}$$ is the dual basis since
clearly
$$(w_{i_1}\wedge \cdots
\wedge w_{i_{\ell}},z_{j_1}\wedge \cdots
\wedge z_{j_{\ell}}) = \prod_{k=1}^n\delta_{i_kj_k}.\eqno (2.6)$$ In addition if the
 $w_i$ are a ${\cal B}$-orthonormal basis of $\g$, then $w_i=z_i,\,i=1,\ldots,n$, and
hence (2.6) implies that $\{w_{i_1}\wedge \cdots
\wedge w_{i_{\ell}}\mid 1\leq i_1< \cdots < i_{\ell}\leq n\}$ is a ${\cal B}$
orthonormal basis of $\wedge^{\ell}\g$. \vs  Now for any
$y_i\in
\g,\,i=1,\ldots,\ell$, let
$\psi(y_1,\ldots,y_{\ell}) = \hbox{\rm det}\,\,
\partial_{y_i}p_j$ so that $$\psi(y_1,\ldots,y_{\ell})\in S^r(\g) \eqno (2.7)$$ by
(2.4). But now $d_Wp_j$ is an invariant 1-form on $\g$. If $x\in \g$, then
$d_Wp_j(x)\in \wedge^1\g$. Explicitly, using the notation in (2.1), one has $$
 d_Wp_j(x) = \sum_{i=1}^n \partial_{z_i}p_j(x)\,w_i.\eqno (2.8)$$ One notes that
$\partial_{z_i}p_j$ is an $n\times \ell$ matrix of polynomial functions. There are
${n\choose \ell}\,\,\ell\times \ell$ minors for this matrix. The determinants of
these minors all lie in $S^r(\g)$ and appear in the following expansion. \vs

 {\bf
Proposition 2.1.} {\it Let the notation be as in (2.1). Let $x\in \g$. Then in
$\wedge^{\ell}\g$ one has $$d_W\,p_1(x)\wedge\cdots \wedge d_W\,p_{\ell}(x) =
\sum_{1\leq i_1< \cdots < i_{\ell}\leq n} \psi(z_{i_1},\ldots,z_{i_{\ell}})(x)
w_{i_1}\wedge \cdots
\wedge w_{i_{\ell}}.\eqno (2.9)$$}\vs {\bf Proof.} This is just standard exterior
algebra calculus using (2.8). QED

\vs {\bf Theorem 2.2.} {\it Let
$v_i,\,i=1,\ldots,n$, be a ${\cal B}$ orthonormal basis of $\g$ chosen and ordered so
that
$v_i,\, i=1,\ldots,\ell$, is a basis of $\hh$. Then there exists a scalar $\kappa \in
\Bbb C^{\times}$ such that, for any $y\in \hh$,
$$d_W\,p_1(y)\wedge\cdots \wedge d_W\,p_{\ell}(y) =
\kappa\,(\prod_{\varphi\in\Delta_+}\varphi(y))\,v_1\wedge\cdots \wedge v_{\ell}.\eqno
(2.10)$$}\vs {\bf Proof.} If $a\in G$, $x\in \g$ and $j = 1,\ldots,\ell$, then since
$d_W\,p_j$ is $G$-invariant one has $$a\cdot d_W\,p_j(x) = d_W\,p_j(a\cdot x).\eqno
(2.11)$$ But this implies that $$d_W\,p_j(x)\in \hbox{\rm cent}\,\g^x\eqno (2.12)$$ since if we
choose
$a\in G^x$ in (2.11) it follows from (2.11) that $d_W\,p_j(x)$ commutes with $\g^x$.
But $x\in \g^x$ so that $d_W\,p_j(x)\in \g^x$. This establishes (2.12).

Now by Theorem 9, p.~382 in [K2] one has that
if $x\in \g$, then $$\{d_W\,p_1(x),\ldots,d_W\,p_{\ell}(x)\}\,\,\hbox{are linearly
independent $\iff\,\,x\in \hbox{\rm Reg}\,\g $}.\eqno (2.12a)$$ Thus the left side of (2.10)
vanishes if and only if
$y\in \hbox{\rm Sing}\,\g\cap \hh$. In particular, choosing the $z_i$ in (2.9)
 so that $v_j= z_j$
for $j=1,\ldots,\ell,$ one has
$\psi(v_1,\ldots, v_{\ell})(y) = 0$ if
$y$ is singular by the expansion (2.9). One the other hand, if $y\in \hh$ is regular
then, by (2.12), one must have that $$\{d_W\,p_j(y),\,j=1,
\ldots,\ell\}\,\,\hbox{is a basis of $\hh$}.\eqno (2.13)$$ Thus if $y$ is regular, the
left side of (2.10) equals $\nu\,v_1\wedge\cdots \wedge v_{\ell}$ for some $\nu\in \Bbb
C^{\times}$. Comparing with the expansion (2.9) one must have
$\nu =
\psi(v_1,\ldots,v_{\ell})(y)$. But then $\psi(v_1,\ldots,v_{\ell})|\hh$ is a polynomial
of of degree $r$ which vanishes on $y\in \hh$ if and only if $y\in \hh$ is singular.
Thus
$$\psi(v_1,\ldots,v_{\ell})|\hh = \kappa \prod_{\varphi\in \Delta_+}\varphi$$ for some
nonzero constant $\kappa$. This proves (2.10). QED

\vs {\bf 2.2.} For any
root $\varphi\in \Delta$ let $e_{\varphi}\in \g$ be a corresponding root vector. We
will make choices so that $$(e_{\varphi}, e_{-\varphi}) = 1.\eqno (2.14)$$ For any $x\in
\hh$, one then has $$dx = \sum_{\varphi\in \Delta_+} \varphi(x)\,e_{\varphi}\wedge
e_{-\varphi}.\eqno (2.15)$$ See Proposition 37, p.~311 in [K5], noting (106), p.~302
and (142), p.~309 in [K5]. But then recalling (1.27) one has $$\gamma_r(x^r) =
r! (-1)^r\,\prod_{\varphi\in \Delta_+} \varphi(x)\,e_{\varphi}\wedge e_{-\varphi}.\eqno
(2.16)$$ But since $(e_{\varphi}\wedge e_{-\varphi},e_{\varphi}\wedge e_{-\varphi}) =
-1$, by (2.14), for any $\varphi\in \Delta_+$ one has that $$(\prod_{\varphi\in
\Delta_+} e_{\varphi}\wedge e_{-\varphi},\prod_{\varphi\in \Delta_+}
\,e_{\varphi}\wedge e_{-\varphi}) = (-1)^r.\eqno (2.17)$$ But then if
$\{v_i\mid \,i=1,\ldots,\ell\}$ is an orthonormal basis of $\hh$, one has
$$(v_1\wedge \cdots\wedge v_{\ell}\wedge \prod_{\varphi\in
\Delta_+} e_{\varphi}\wedge e_{-\varphi},v_1\wedge \cdots\wedge v_{\ell}\wedge
\prod_{\varphi\in \Delta_+}
\,e_{\varphi}\wedge e_{-\varphi}) = (-1)^r. \eqno (2.18)$$ But then we may choose an
ordering of the $v_i$ such that $$\mu = i^r v_1\wedge \cdots\wedge v_{\ell}\wedge
\prod_{\varphi\in
\Delta_+} e_{\varphi}\wedge e_{-\varphi}\eqno (2.19)$$ so that $$(v_1\wedge\cdots
\wedge v_{\ell})^* = i^r \prod_{\varphi\in
\Delta_+} e_{\varphi}\wedge e_{-\varphi}.\eqno (2.20)$$ But then one has

\vs {\bf Theorem 2.3} {\it There exists $\kappa_o\in \Bbb C^{\times}$ such that for any
$x\in \g$, $$\eqalign{(d_W\,p_1(x)\wedge\cdots \wedge d_W\,p_{\ell}(x))^* &=
\kappa_o\,\,{(-dx)^r\over r!}\cr &=\kappa_o\,\, \gamma_r({x^r\over r!}).\cr}\eqno
(2.21)$$}\vs {\bf Proof.} If $y\in\hh$ is regular, then (2.21), for $y=x$, follows from
(2.16),(2.20) and Theorem 2.2. That is $$\eqalign{(d_W\,p_1(y)\wedge\cdots \wedge
d_W\,p_{\ell}(y))^* &=
\kappa_o\,\,{(-dy)^r\over r!}\cr &=\kappa_o\,\, \gamma_r({y^r\over r!}).\cr}\eqno
(2.22)$$ But now if $x\in \g$ regular and semisimple there exist $a\in G$ and a regular
$y\in \hh$ such that $a\cdot y = x$. But now since $^*$ and $\gamma_r$ are clearly
$G$-maps one has (2.21) by applying the action of $a$ to both sides of (2.22). However
the set of regular semisimple elements in $\g$ is dense (this nonempty set is Zariski
open) one has (2.21) for all $x\in \g$ by continuity. QED

\vs  Returning to
our module $M$ of harmonic polynomials on $\g$ of degree $r$ it is obvious, by
definition, that $M$ is spanned by all $f\in S^r$ of the form $f = \Gamma (w_1\wedge \cdots
\wedge  w_{2r})$ where the $w_i\in \g$ are linearly independent. Explicitly
$\Gamma (w_1\wedge \cdots \wedge w_{2r})$ is given by (1.22). We now show that $\Gamma
(w_1\wedge \cdots \wedge w_{2r})$ may also be given as the determinant of one of the
$\ell\times \ell$ minors in the expansion (2.9).

\vs {\bf Theorem 2.4.} {\it Let
$w_k\in \g, \, k=1,\ldots, 2r,$ be linearly independent and let $\ss\s\g$ be the
span of the $w_k$ and let $u_i\in \g, =1,\ldots, \ell$, be a basis of $\ss^{\perp}$.
Then there exists a constant $\kappa_1\in \Bbb C^{\times}$ such that $$\eqalign{\Gamma
(w_1\wedge
\cdots \wedge w_{2r}) &= \kappa_1\,\psi(u_1,\ldots, u_{\ell})\cr &=
\kappa_1\, \hbox{\rm det}\,\partial_{u_i}\,p_j.\cr}\eqno (2.23)$$ \vskip .5pc Furthermore $M$
is the span of all $\ell\times \ell$ determinant minors $\psi(v_1,\ldots,
v_{\ell})$ where $v_i\in\g,\,i=1,\ldots,\ell$, are linearly independent.}

\vs {\bf
Proof.} Clearly we may choose the two dual bases in (2.1) so that the given
$w_k$ are the first
$2r$-elements of the $w$ basis and the $u_i$ are the last $\ell$ elements of the $z$
basis. Thus there exists $\kappa_2\in \Bbb C^{\times} $ such that
$$(u_1\wedge\cdots\wedge u_{\ell})^* = \kappa_2\,w_1\wedge \cdots \wedge w_{2r}.\eqno
(2.24)$$ Now let $x\in \g$. Then by the expansion (2.9) one has
$$(d_W\,p_1(x)\wedge\cdots \wedge d_W\,p_{\ell}(x), u_1\wedge \cdots\wedge u_{\ell})=
\psi(u_1,\ldots,u_{\ell})(x).\eqno (2.25)$$ But then by (1.3), (1.14), (1.29) and (2.21) one
has
$$\eqalign{\psi(u_1,\ldots,u_{\ell})(x)&= ((d_W\,p_1(x)\wedge\cdots \wedge
d_W\,p_{\ell}(x))^*, (u_1\wedge \cdots\wedge u_{\ell})^*)\cr
&=\kappa_o\kappa_2\,(\gamma_r({x^r\over r!}), w_1\wedge \cdots \wedge w_{2r})\cr
&= \kappa_1^{-1}\,\Gamma(w_1\wedge \cdots \wedge w_{2r})(x)\cr}\eqno (2.26)$$ where
$\kappa_1^{-1} = \kappa_o\kappa_2$. The last statement in the theorem is obvious since
clearly $u_i,\, i=1,\ldots, \ell,$ is an aritrary set of $\ell$-independent elements
in $\g$.  QED

\vs {\bf 2.3.} Let $\{z_i,w_j\}$ be the arbitrary dual bases of $\g$ as
in (1.4). Then, independent of the choice of bases, the Casimir operator $\hbox{\rm Cas}$ on
$\wedge^\g$ is given by $$\hbox{\rm Cas} = \sum_{i=1}^n \theta(z_i)\theta(w_i).$$ We recall
special cases of some results in [K3]. Let $A_{\ell}\s \wedge^{\ell}\g$ be the span
in $\wedge^{\ell}\g$ of all $[\cc]$ where $\cc\in \g$ is a commutative Lie subalgebra of
dimension $\ell$. Since the set of such subalgebras includes, for example, Cartan
subalgebras it is obvious that $A_{\ell}\neq 0$. In fact note that $$[\g^y]\s
A_{\ell}\eqno (2.27)$$ for any $y\in \hbox{\rm Reg}\,\g$ since, as one knows, $\g^y$ is abelian
if $y$ is regular. Clearly
$A_{\ell}$ is a $G$ submodule of $\wedge^{\ell}$. On the other hand, let
$m_{\ell}$ be the maximal value of $\hbox{\rm Cas}$ on $\wedge^{\ell}$ and let $M_{\ell}$ be the
corresponding $\hbox{\rm Cas}$ eigenspace. Again, clearly $M_{\ell}$ is a  $G$-submodule
of
$\wedge^{\ell}\g$. From the definition of $M_{\ell}$ it is obvious that
$\hbox{\rm Hom}_G(M_{\ell}, \wedge^{\ell}\g/M_{\ell}) = 0$. Since ${\cal B}|\wedge^{\ell}\g$ is
nonsingular it follows that $${\cal B}|M_{\ell}\,\,\hbox{is nonsingular}\eqno
(2.28)$$ and hence $M_{\ell}$ is self-contragredient.
Noting the $1/2$ in (2.1.7) of [K3] the following result is a special case of Theorem (5), p.
156 in [K3].

 \vs {\bf Theorem 2.5.} {\it One has
$$A_{\ell} = M_{\ell} \eqno (2.29)$$ and in addition $$m_{\ell} = \ell.\eqno (2.30)$$}\vs For any
ordered subset $\Phi\s \Delta$, $\Phi =\{\varphi_1,\ldots,\varphi_k\}$, let $e_{\Phi} =
e_{\varphi_1}\wedge \cdots\wedge e_{\varphi_k}$ and put $\langle \Phi\rangle = \sum_{\varphi\in
\Phi} \varphi $ so that with respect to $\hh$, $$e_{\Phi}\in \wedge^k\g\,\,\,\hbox{is a weight
vector of weight $\langle \Phi\rangle$}.\eqno (2.31)$$\vskip .5pc Let $\b\s\g$ be the Borel
subalgebra of $\g$ spanned by $\hh$ and $\{e_{\varphi}\}$, for $\varphi\in \Delta_+$, and put
 $\n = [\b,\b]$. Any ideal $\a$ of $\b$ where $\a\s \n$ is necessarily spanned by root vectors. We
will say that $\Phi$, as above, is an ideal of $\Delta_+$ if $\Phi\s \Delta_+$ and $\a_{\Phi} =
\sum_{i=1}^{k} \Bbb C e_{\varphi_i}$ is an ideal in $\b$.

\vs {\bf Remark 2.6.} One notes that if
$\Phi$ is an ideal of $\Delta_+$ and $V_{\Phi}\s \wedge^k\g$ is the  $G$-module spanned
by $G\cdot e_{\Phi}$, then $V_{\Phi}$ is irreducible having $e_{\Phi}$ as highest weight vector
and $\langle \Phi\rangle$ as highest weight.

As already noted in [K3] (see bottom of p. 158) it is immediate that if $\a$ is any abelian ideal
in $\b$, then $\a \s \n$ so that $\a = \a_{\Phi}$ for an ideal $\Phi\s \Delta_+$. Much more
subtly it has been established in [KW] (see Lemma 12, p. 113 in [KW]) that any ideal $\a$ of $\b$
having dimension $\ell$ is in fact abelian. Let ${\cal I}$ be the (obviously finite) set
 of  all ideals
$\Phi$ in $\Delta_+$ which have cardinality $\ell$. If $\Phi_1,\Phi_2\in {\cal I}$ are distinct,
then $\langle \Phi_1\rangle  \neq \langle \Phi_2\rangle$ by Theorem (7), p.~158 in [K3] so that
$V_{\Phi_1}$ are inequivalent $\g$ and $G$ modules. Then Theorem (8), p. 159 in [K3] implies \vs

{\bf Theorem 2.7.} {\it $M_{\ell}$ is a multiplicity one  $G$-module. In fact $$M_{\ell}
=
\oplus_{\Phi\in {\cal I}}\,\, V_{\Phi}\eqno (2.32)$$ so the number of irreducible components in
$M_{\ell}$ is the cardinality of ${\cal I}$. }\vs

{\bf Remark 2.8.} In the general case
we do not have a formula for $\hbox{\rm card}\,{\cal I}$ although computing this number in any
given case does not seem to be too difficult. In the special case where $\g\cong
\hbox{\rm Lie} \,Sl(n,\Bbb C)$ one easily has a bijective correspondence of ${\cal I}$ with
the  set of
all Young tableaux of size $n-1$ so that in this case $$\hbox{\rm card}\,{\cal I} = p(n-1)\eqno
(2.33)$$ where $p$ here is the classical partition  function.\vs
 Let
$$\tau:\wedge^{\ell}\g\to
\wedge^{2r}\g\eqno (2.34)$$ be the $G$-ismorphism defined by putting $\tau(u) = u^*$
recalling that $u^*=
\iota(u)\mu$. Let $M_{2r} = \tau(M_{\ell})$, \vs

{\bf Theorem 2.9.} {\it $\tau$ is a ${\cal
B}$-isomorphism so that ${\cal B}|M_{2r}$ is nonsingular. Furthermore $\ell$ is the maximal
eigenvalue of $\hbox{\rm Cas}$ on $\wedge^{2r}\g$ and $M_{2r}$ is the corresponding  eigenspace.
As
$G$ modules one has $$M_{\ell}\cong M_{2r}\eqno (2.35)$$ so that $M_{2r}$ is a
multiplicity 1 module where in fact $$M_{2r}\cong
\oplus_{\Phi\in {\cal I}}\,\, V_{\Phi}.\eqno (2.36)$$
 We recall the $V_{\Phi}$ is an irreducible
$G$-module with highest weight $\langle \Phi \rangle$. See (2.31).}\vs {\bf Proof.} The first
statement follows from Proposition 1.1. The remaining statements are immediate from Theorem 2.7
since $\tau$ is a $G$-isomorphism. QED

\vs In  light of equality $M_{\ell}= A_{\ell}$ (see
(2.9)) Ranee Brylinski in her thesis (see [RB]) proved that $M_{\ell}$ is the span of $G\cdot
[\hh]$. The thesis however has not been published. A stronger theorem (motivated by her result)
appears in [KW]. The following result is just Corollary 2, p. 105 in [KW].

\vs {\bf Theorem 2.10.}
{\it $M_{\ell}$ is the span of $G\cdot [\g^x]$ for any $x\in \hbox{\rm Reg}\,\g$.}\vs Now by
(2.12) and (2.12a) one has $$\Bbb C\,\,d_W\,p_1(x)\wedge\cdots \wedge d_W\,p_{\ell}(x) =
[\g^x]\eqno (2.37)$$ for any
$x\in \hbox{\rm Reg}\,\g$. Using Theorem 2.3 we can now transfer Theorem 2.10 to $M_{2r}$ where it
will have consequences for the structure of the space of functions $M\s H^r$.

\vs {\bf Theorem
2.11.} {\it
$M_{2r}$ is the span of $G\cdot (\gamma_r({x^r \over r!})) $ for any $x\in \hbox{\rm
Reg}\,\g$.}\vs Proof. This is immediate from Theorem 2.3, Theorem 2.10, (2.37) and the fact that
$\tau$ is a
$G$-isomorphism. QED. \vs Let $N_{2r}$ be the ${\cal B}$ orthogonal subspace to $M_{2r}$ in
$\wedge^{2r}\g$. By the first statement in Theorem 2.9 one has a ${\cal B}$ orthogonal $G$-module
decomposition
$\wedge^{2r}\g$, $$\wedge^{2r}\g = N_{2r}\oplus M_{2r}.\eqno (2.38)$$

\vskip .5pc {\bf Remark 2.12.} Note that by Theorem 2.9 any eigenvalue of $\hbox{\rm Cas}$ in
$N_{2r}$ is less than $\ell$. \vs We return now to our $G$-space $M$ of homogeneous harmonic
polynomials on
$\g$ of degree $r$ which define $\hbox{\rm Sing}\,\g$. We recapitulate some of the properties of
$M=
\Gamma(\wedge^{2r}\g)$ already established in this paper. Let $w_k\in \g,\,k=1,\ldots,2r$, be
linearly independent and let $z_i\in \g,\,\,i=1,\ldots,\ell,$ be linearly independent and ${\cal B}$
orthogonal to the $w_k$. Then for suitable generators $p_j,\,j=1,\ldots,\ell,$ of $J= S(\g)^G$,
we have
$$\eqalign{(1)&\,\,
\Gamma(w_i\wedge\cdots\wedge w_{2r})\,\,\hbox{is explicitly given by (1.23)}\cr (2)&\,\,
\Gamma(w_i\wedge\cdots\wedge w_{2r})\,\hbox{is given as (up to scalar multiplication)
$det\,\partial_{z_i}p_j$}.\,\, \hbox{See Theorem 2.4.}\cr (3)&\,\,\hbox{If $f\in
M$, then $f|\a$, where
$\a$ is any Cartan subalgebra or}\cr &\,\,\,\,\,\,\,
\hbox{$\a=\g^e$ for $e$ principal nilpotent, is given in Theorems 1.6 and 1.7.}\cr}$$\vskip .5pc
We now determine the $G$-module structure of $M$, \vs {\bf Theorem 2.13.} {\it
$N_{2r} =
\hbox{\rm Ker}\,\Gamma$ and
$$\Gamma:M_{2r} \to M\eqno (2.39)$$ is a $G$-isomorphism so that as $G$-modules $$M \cong
M_{2r}\cong M_{\ell} = A_{\ell}\eqno (2.40)$$ where we recall $A_{\ell}\s
\wedge^{\ell}\g$ has been defined in [K3] as the span of $[\ss]$ over all abelian
subalgebras $\ss\s \g$ of dimension $\ell$.

Furthermore we have defined ${\cal I}$ as the set of all ideals $\Phi$ in $\Delta_+$ of
cardinality $\ell$, parameterizing with the notation $\a_{\Phi}$, the set of all ideals $\a$ of
$\b$ having dimension $\ell$. See Remark 2.6.

Moreover $M$ is a multiplicity one $G$-module with $\hbox{\rm card}\,{\cal I}$ irreducible
components. In addition ${\cal I}$ parameterizes these components in the sense that the component
corresponding to
$\Phi\in {\cal I}$ is equivalent to $V_{\Phi}$, using the notation of Remark 2.6, and hence has
highest weight $\langle\Phi\rangle$. Finally $Cas$ takes the value $\ell$ on each and every
irreducible component of $M$.}

\vs {\bf Proof.} By (1.27) and (1.29) one has $$(\Gamma(\zeta)(x) =
(\zeta,\gamma_r({x^r \over r!}))\eqno (2.41)$$ for any $x\in \g$ and any $\zeta\in\wedge^{2r}\g$.
Of course $\gamma_r({x^r \over r!})= 0$ for any $x\in \hbox{\rm Sing}\,\g$ (see (2.12a) and
Theorem 2.3). However
$M_{2r}$ is the span of $G\cdot \gamma_r({x^r \over r!})$ for any $x\in \hbox{\rm Reg}\,\g$ by
Theorem 2.11. Thus not only does (2.41) imply that $N_{2r}\s \hbox{\rm Ker}\,\Gamma $ but $N_{2r}
= \hbox{\rm Ker} \,\Gamma$ since if
$\zeta\in M_{2r}$ and $x\in \hbox{\rm Reg}\,\g $ there exists $a\in G$ such that if $y=a\cdot x$,
then
$\Gamma(\zeta)(y)\neq 0$ by Theorem 2.11 and the nonsingularity of ${\cal B}|M_{2r}$, as asserted
in  Theorem 2.9. Since $\Gamma$ is a $G$-map one has the isomorphism (2.39). The remaining
statements follow from Theorem 2.5 and Theorem 2.9. QED

\centerline{\bf References}\vskip .5pc

\item {[K1]} B. Kostant,  The Three Dimensional Sub-Group and the Betti Numbers of a Complex
Simple Lie  Group, {\it Amer. Jour. of Math.}, {\bf 81}(1959), 973--1032.
\item {[K2]} B. Kostant,  Lie Group Representations on Polynomial Rings, {\it Amer. J. Math.},
  {\bf 85}(1963), 327--404.
\item {[K3]} B. Kostant,  Eigenvalues  of a Laplacian and Commutative Lie Subalgebras,
 {\it Topology}, {\bf 13}(1965) 147--159.
\item {[K4]} B. Kostant,  A Lie Algebra Generalization of the Amitsur-Levitski Theorem, {\it Adv. in
Math.}, {\bf 40}(1981), No. 2, 155--175.
\item {[K5]} B. Kostant,  Clifford Algebra Analogue of the
Hopf-Koszul-Samelson Theorem, the $\rho$-Decomposition,
$C(\g) = \hbox{\rm End}\,V_{\rho}\otimes C(P)$, and the $\g$-Module Structure of $\wedge \g$, {\it Adv.
in Math.}, {\bf 125}(1997), 275--350.
\item {[KW]}   B. Kostant and N. Wallach,  On a theorem of Ranee Brylinski,
 {\it Contemporary Mathematics}, {\bf 490}(2009),
105--142.
\item {[Kz]} J. L. Koszul, Homologie et cohomologie des alg\`ebres de Lie, {\it  Bull. Soc.
Math. Fr.}, {\bf 78}(1950), 65--127.
\item {[RB]} R. Brylinski, {\it Abelian algebras and adjoint orbits}, Thesis MIT, 1981.

\end

\end